\documentstyle[11pt]{article}

\input amssym.def
\input amssym.tex

\newtheorem{theorem}{Theorem}[section] 
\newtheorem{question}[theorem]{Problem} 
\newtheorem{definition}[theorem]{Definition}
\newtheorem{thesis}[theorem]{Thesis} 
\newtheorem{remark}[theorem]{Remark}

\newcommand{\cf}{{\rm cf}}
\newcommand{\pp}{{\rm pp}}
\newcommand{\Levy}{{\rm Levy}}
\newcommand{\bd}{{\rm bd}}
\newcommand{\Ga}{{\frak a}}
\newcommand{\Reg}{{\rm Reg}}
\newcommand{\pcf}{{\rm pcf}}
\newcommand{\cov}{{\rm cov}}
\newcommand{\acc}{{\rm acc}}
\newcommand{\otp}{{\rm otp}}
\newcommand{\Proof}{\noindent {\sc Proof}\hskip0.2in}
\newcommand{\nacc}{{\rm nacc}}
\newcommand{\qed}{\hspace{0.2in}\vrule width 6pt height 6pt depth 0pt 
\vspace{0.1in}} 

\title{A polarized partition relation and failure of GCH at singular strong
limit}
\author{{\bf Saharon Shelah}\thanks{\null\ Research partially supported by
``Basic Research Foundation'' administered by The Israel Academy of Sciences
and Humanities. Publication 586.}\\
Institute of Mathematics\\
The Hebrew University\\
91 904 Jerusalem, Israel\\
and\\
Department of Mathematics\\
Rutgers University\\
New Brunswick, N.J. 08903, USA}
\date{\today}

\begin{document}
\baselineskip13.14 truept
\maketitle
\begin{abstract}
The main result is that for $\lambda$ strong limit singular failing
the continuum hypothesis (i.e. $2^\lambda> \lambda^+$), a polarized
partition theorem holds.
\end{abstract}

\section{Introduction}
In the present paper we show a polarized partition theorem for strong limit
singular cardinals $\lambda$ failing the continuum hypothesis. Let us recall
the following definition. 
\begin{definition}
\label{jeden}
For ordinal numbers $\alpha_1,\alpha_2,\beta_1,\beta_2$ and a cardinal
$\theta$, the polarized partition symbol
$$
\left( \begin{array}{c} \alpha_1 \\ \beta_1 \end{array}
\right) \rightarrow \left( \begin{array}{c} \alpha_2 \\ \beta_2
\end{array} \right)^{1,1}_{\theta}
$$
means:\\
{\em if} $d$ is a function from $\alpha_1\times \beta_1$ into $\theta$
{\em then} for some $A\subseteq \alpha_1$ of order type $\alpha_2$ and
$B\subseteq \beta_1$ of order type $\beta_2$, the function $d\restriction
A\times B$ is constant. 
\end{definition}
We address the following problem of Erd\H os and Hajnal:
\begin{description}
\item[(*)] if $\mu$ is strong limit singular of uncountable cofinality, 
$\theta<\cf(\mu)$ does
\[\left( \begin{array}{c}
\mu^+\\ \mu 
\end{array} \right) \rightarrow \left( \begin{array}{c}
\mu \\ \mu \end{array} \right)^{1,1}_\theta\qquad ?\] 
\end{description}
The particular case of this question for $\mu=\aleph_{\omega_1}$ and
$\theta=2$ was posed by Erd\H os, Hajnal and Rado (under the assumption of
GCH) in \cite[Problem 11, p.183]{EHR}). Hajnal said that the assumption of GCH
in \cite{EHR} was not crucial, and he added that the intention was to ask the
question ``in some, preferably nice, Set Theory''. 

Baumgartner and Hajnal have proved that if $\mu$ is weakly compact then the
answer to {\bf (*)} is ``yes'' (see \cite{BH}), also if $\mu$ is strong limit
of cofinality $\aleph_0$. But for a weakly compact $\mu$ we do not know
if for every $\alpha<\mu^+$: 
\[\left( \begin{array}{c} \mu^+ \\
\mu \end{array} \right) \rightarrow \left( \begin{array}{c} \alpha \\
\mu \end{array} \right)^{1,1}_{\theta}.\]
The first time I heard the problem (around 1990) I noted that {\bf (*)} holds
when $\mu$ is a singular limit of measurable cardinals. This result is
presented in Theorem \ref{piec}. It seemed likely that we can combine this
with suitable collapses, to get ``small" such $\mu$ (like $\aleph_{\omega_1}$)
but there was no success in this direction. 

In September 1994, Hajnal reasked me the question putting great stress
on it. Here we answer the problem {\bf (*)} using methods of \cite{Sh:g}. 
But instead of the assumption of GCH (postulated in \cite{EHR}) we assume
$2^\mu>\mu^+$. The proof seems quite flexible but we did not find out what
else it is good for. This is a good example of the major theme of \cite{Sh:g}:

\begin{thesis}
Whereas CH and GCH are good (helpful, strategic) assumptions having many
consequences, and, say, $\neg$CH is not, the negation of GCH at singular
cardinals (i.e. for $\mu$ strong limit singular $2^\mu>\mu^+$ or, really the
strong hypothesis: $\cf(\mu)<\mu\quad\Rightarrow\quad\pp(\mu)>\mu^+$) is a
good (helpful, strategic) assumption.
\end{thesis}
Foreman pointed out that the result presented in Theorem \ref{cztery} below is
preserved by $\mu^+$-closed forcing notions. Therefore, if   
\[V \models 
\left( \begin{array}{c} \lambda^+ \\ \lambda \end{array}\right)
\rightarrow \left( \begin{array}{c} \lambda \\ \lambda \end{array}
\right)^{1,1}_\theta\]
then
\[V^{\Levy (\lambda^+, 2^\lambda)} \models 
\left( \begin{array}{c} \lambda^+ \\ \lambda \end{array} \right)
\rightarrow \left(\begin{array}{c} \lambda \\ \lambda \end{array}
\right)^{1,1}_\theta .\]
Consequently, the result is consistent with $2^\lambda = \lambda^+ \ \&\
\lambda$ is small. (Note that although our final model may satisfy the
Singular Cardinals Hypothesis, the intermediate model still violates SCH at
$\lambda$, hence needs large cardinals, see \cite{J}.) For $\lambda$ not small
we can use Theorem \ref{piec}). 

Before we move to the main theorem, let us recall an open problem important
for our methods: 
\begin{question}
\label{trzy}
\begin{enumerate}
\item Let $\kappa=\cf(\mu)>\aleph_0$, $\lambda$, $\mu>2^\theta$ and
$\lambda=\cf(\lambda)\in (\mu, \pp^+(\mu))$. Can we find
$\theta<\mu$ and $\Ga\in [\mu\cap \Reg]^\theta$ such that: $\lambda\in
\pcf(\Ga)$, $\Ga=\bigcup\limits_{i<\kappa}\Ga_i$, $\Ga_i$ bounded in
$\mu$ and
$\sigma\in \Ga_i\Rightarrow \bigwedge\limits_{\alpha<\sigma}
|\alpha|^{\theta}< \sigma$?

For this it is enough to show:
\item If $\mu=\cf(\mu)> 2^{<\theta}$ but $\bigvee\limits_{\alpha<\mu}
|\alpha|^{<\theta}\geq \mu$ then we can find
$\Ga\in [\mu\cap \Reg]^{<\theta}$ such that $\lambda\in \pcf(\Ga)$.
\end{enumerate}
\end{question}
As shown in \cite{Sh:g}
\begin{theorem}
\label{dwa}
If $\mu$ is strong limit singular of cofinality $\kappa>\aleph_0$,
$2^\mu > \lambda=\cf(\lambda)> \mu$ then for some strictly increasing sequence
$\langle \lambda_i:i<\kappa\rangle$ of regulars with limit $\mu$,
$\prod\limits_{i<\kappa}\lambda_i/J^{\bd}_\kappa$ has true cofinality
$\lambda$. If $\kappa=\aleph_0$, it still holds for $\lambda=\mu^{++}$.
\end{theorem}
[More fully, by \cite[II \S5]{Sh:g}, we know $\pp(\mu)=^+ 2^\mu$ and by
\cite[VIII 1.6(2)]{Sh:g}, we know $\pp^+(\mu)=\pp^+_{J^{\bd}_\kappa}(\mu)$.
Note that for $\kappa=\aleph_0$ we should replace $J^{\bd}_\kappa$ by a
possibly larger ideal, using \cite[1.1, 6.5]{Sh:430} but there is no need
here.]  

\begin{remark}
{\em
Note the problem is $\pp=\cov$ problem, see more \cite[\S1]{Sh:430}; so
if $\kappa=\aleph_0$, $\lambda<\mu^{+\omega_1}$ the conclusion of
\ref{dwa} holds; we allow to increase $J^{\bd}_\kappa$, even ``there
are $<\mu^+$ fixed points $<\lambda^+$'' suffices.
}
\end{remark}

\section{Main result}

\begin{theorem}
\label{cztery}
Suppose $\mu$ is strong limit singular satisfying $2^\mu>\mu^+$. Then
\begin{enumerate}
\item $\left( \begin{array}{c} \mu^+ \\ \mu\end{array} \right)
\rightarrow \left(\begin{array}{c} \mu+1 \\ \mu \end{array}
\right)^{1,1}_{\theta}$ for any $\theta< \cf(\mu)$,
\item if $d$ is a function from $\mu^+\times\mu$ to $\theta$ and
$\theta<\cf(\mu)$ then for some sets $A\subseteq\mu^+$ and $B\subseteq\mu$ we
have:\qquad $\otp (A)=\mu+1$, $\otp (B)=\mu$ and the restriction
$d\restriction A\times B$ does not depend on the first coordinate. 
\end{enumerate}
\end{theorem}

\Proof 1) It follows from part 2) (as if $d(\alpha,\beta)= d'(\beta)$ for
$\alpha\in A$, $\beta\in B$, where $d': B\rightarrow\theta$, and $|B|=\mu$,
$\theta<\cf(\mu)$ then there is $B'\subseteq B$, $|B'|=\mu$ such that
$d'\restriction B$ is constant and hence $d\restriction (A\times B')$ is
constant as required).

\noindent 2) Let $d:\mu^+\times\mu\rightarrow\theta$. Let $\kappa=\cf(\mu)$
and $\bar{\mu}=\langle\mu_i:i<\kappa \rangle$ be a continuous strictly
increasing sequence such that $\mu=\sum\limits_{i<\kappa} \mu_i$,  $\mu_0>
\kappa$. We can find a sequence $\bar C=\langle C_\alpha:\alpha<\mu^+\rangle$
such that: 
\begin{description}
\item[(A)] $C_\alpha\subseteq \alpha$ is closed, $\otp (C_\alpha)<\mu$,
\item[(B)] $\beta\in \nacc (C_\alpha)\quad \Rightarrow\quad C_\beta=C_\alpha
\cap\beta$,
\item[(C)] if $C_\alpha$ has no last element then $\alpha=\sup(C_\alpha)$ (so
$\alpha$ is a limit ordinal) and any member of $\nacc(C_\alpha)$ is a
successor ordinal,
\item[(D)] if $\sigma=\cf(\sigma)<\mu$ then the set 
\[S_\sigma=:\{\delta<\mu^+:\cf(\delta)=\sigma\ \ \&\ \ \delta=\sup(C_\delta)\
\ \&\ \ \otp(C_\delta)=\sigma\}\]
is stationary
\end{description}
(possible by \cite[\S1]{Sh:420}); we could have added 
\begin{description}
\item[(E)] for every $\sigma\in\Reg\cap \mu^+$ and a club $E$ of $\mu^+$, for
stationary many $\delta\in S_\sigma$,  $E$ separates any two successive
members of $C_\delta$.
\end{description}
Let $c$ be a symmetric two place function from $\mu^+$ to $\kappa$ such that
for each $i<\kappa$ and $\beta<\mu^+$ the set
\[a^\beta_i=: \{\alpha<\beta: c(\alpha, \beta)\leq i\}\]
has cardinality $\leq \mu_i$ and $\alpha<\beta<\gamma \Rightarrow c(\alpha,
\gamma)\leq \max\{ c(\alpha, \beta), c(\beta, \gamma)\}$ and
\[\alpha\in C_\beta \ \&\ \mu_i \geq |C_\beta|\quad\Rightarrow\quad c(\alpha,
\beta)\leq i\]
(as in \cite{Sh:108}, easily constructed by induction on $\beta$).

Let $\bar \lambda=\langle \lambda_i: i<\kappa\rangle$ be a strictly increasing
sequence of regular cardinals with limit $\mu$ such that
$\prod\limits_{i<\kappa}\lambda_i/ J^{\bd}_\kappa$ has true cofinality
$\mu^{++}$ (exists by \ref{dwa} with $\lambda=\mu^{++}\leq 2^\mu$). As we can
replace $\bar \lambda$ by any subsequence of length $\kappa$, without loss of
generality $(\forall i<\kappa)(\lambda_i> 2^{\mu^+_i})$. Lastly, let $\chi=
\beth_8(\mu)^+$ and $<^*_\chi$ be a well ordering of ${\cal H}(\chi)$($=:\{x$:
the transitive closure of $x$ is of cardinality $<\chi\}$).

Now we choose by induction on $\alpha<\mu^+$ sequences $\bar M_\alpha=\langle
M_{\alpha, i}: i< \kappa\rangle$ such that: 
\begin{description}
\item[(i)]   $M_{\alpha,i}\prec (H(\chi),\in,<^*_\chi)$,
\item[(ii)]  $\|M_{\alpha,i}\|=2^{\mu^+_i}$ and ${}^{\mu^+_i}(M_{\alpha,i})
\subseteq M_{\alpha,i}$ and $2^{\mu^+_i}+1\subseteq M_{\alpha,i}$,
\item[(iii)] $d,c,\bar C,\bar\lambda,\bar\mu,\alpha\in M_{\alpha,i}$,
$\langle M_{\beta,j}:\beta<\alpha,j<\kappa\rangle\in M_{\alpha,i}$,
$\bigcup\limits_{\beta\in a^\alpha_i} M_{\beta,i}\subseteq M_{\alpha,i}$ and
$\langle M_{\alpha,j}: j<i\rangle\in M_{\alpha,i}$, $\bigcup\limits_{j<i}
M_{\alpha,j}\subseteq M_{\alpha,i}$, 
\item[(iv)]  $\langle M_{\beta,i}:\beta\in a^\alpha_i\rangle$ belongs to
$M_{\alpha, i}$. 
\end{description}
There is no problem to carry out the construction. Note that actually the
clause (iv) follows from (i)--(iii), as $a^\alpha_i$ is defined from $c$,
$\alpha$, $i$. Our demands imply that   
\[[\beta\in a^\alpha_i\quad\Rightarrow\quad M_{\beta,i}\prec M_{\alpha,i}]
\qquad\mbox{and}\qquad[j<i\quad\Rightarrow\quad M_{\alpha,j}\prec M_{\alpha,
i}]\] 
and $a^\alpha_i\subseteq M_{\alpha,i}$, hence $\alpha\subseteq
\bigcup\limits_{i<\kappa} M_{\alpha,i}$. 

For $\alpha< \mu^+$ let $f_\alpha\in\prod\limits_{i<\kappa}\lambda_i$ be
defined by $f_\alpha(i)=\sup(\lambda_i\cap M_{\alpha,i})$. Note that
$f_\alpha(i)<\lambda_i$ as $\lambda_i=\cf(\lambda_i)> 2^{\mu^+_i}=\|M_{\alpha,
i}\|$. Also, if $\beta<\alpha$ then for every $i\in [c(\beta,\alpha),\kappa)$
we have $\beta\in M_{\alpha, i}$ and hence $\bar M_\beta\in M_{\alpha,
i}$. Therefore, as also $\bar\lambda\in M_{\alpha, i}$, we have $f_\beta\in
M_{\alpha, i}$ and $f_\beta(i)\in M_{\alpha, i}\cap \lambda_i$. Consequently
\[(\forall i\in[c(\beta,\alpha),\kappa))(f_\beta(i)<f_\alpha(i))\qquad\mbox{
and thus}\quad f_\beta<_{J^{\bd}_\kappa} f_\alpha.\]
Since $\{f_\alpha: \alpha<\mu^+\}\subseteq \prod\limits_{i<\kappa}\lambda_i$
has cardinality $\mu^+$ and $\prod\limits_{i<\kappa}\lambda_i /J^{\bd}_\kappa$
is $\mu^{++}$-directed, there is $f^*\in \prod\limits_{i<\kappa} \lambda_i$
such that 
\begin{description}
\item[$(*)_1$] $(\forall\alpha<\mu^+)(f_\alpha<_{J^{\bd}_\kappa}f^*)$.
\end{description}
Let, for $\alpha<\mu^+$, $g_\alpha\in {}^\kappa \theta$ be defined by 
$g_\alpha(i)=d(\alpha, f^*(i))$. Since $|{}^\kappa\theta|<\mu<\mu^+=
\cf(\mu^+)$, there is a function $g^*\in {}^\kappa\theta$ such that
\begin{description}
\item[$(*)_2$] the set $A^*=\{\alpha<\mu^+: g_\alpha=g^*\}$ is unbounded in
$\mu^+$.
\end{description}
Now choose, by induction on $\zeta<\mu^+$, models $N_\zeta$ such that:
\begin{description}
\item[(a)] $N_\zeta\prec (H(\chi), \in, <^*_\chi)$,
\item[(b)] the sequence $\langle N_\zeta:\zeta<\mu^+\rangle$ is increasing
continuous,
\item[(c)] $\|N_\zeta\|=\mu$ and ${}^{\kappa>}(N_\zeta)\subseteq N_\zeta$
if $\zeta$ is not a limit ordinal,
\item[(d)] $\langle N_\xi:\xi\leq \zeta\rangle \in N_{\zeta+1}$,
\item[(e)] $\mu+1\subseteq N_\zeta$, $\bigcup\limits_{\scriptstyle \alpha<
\zeta\atop\scriptstyle i<\kappa} M_{\alpha,i}\subseteq N_\zeta$ and $\langle
M_{\alpha, i}:\alpha<\mu^+, i<\kappa\rangle$, $\langle f_\alpha: \alpha<\mu^+
\rangle$, $g^*$, $A^*$ and $d$ belong to the first model $N_0$.
\end{description}
Let $E=:\{\zeta<\mu^+: N_\zeta\cap \mu^+=\zeta\}$. Clearly, $E$ is a club of
$\mu^+$, and thus we can find an increasing sequence $\langle\delta_i:
i<\kappa\rangle$ such that 
\begin{description}
\item[$(*)_3$] $\delta_i\in S_{\mu^+_i}\cap \acc (E)$ ($\subseteq\mu^+$) (see
clause (D) in the beginning of the proof).
\end{description}
For each $i<\kappa$ choose a successor ordinal $\alpha^*_i\in\nacc(
C_{\delta_i})\setminus\bigcup \{\delta_j+1:j<i\}$. Take any $\alpha^*\in A^*
\setminus \bigcup\limits_{i<\kappa} \delta_i$.

\noindent We choose by induction on $i<\kappa$ an ordinal $j_i$ and sets
$A_i$, $B_i$ such that: 
\begin{description}
\item[$(\alpha)$] $j_i<\kappa$ and $\mu_{j_i}>\lambda_i$ (so $j_i>i$) and
$j_i$ strictly increasing in $i$,
\item[$(\beta)$] $f_{\delta_i}\restriction [j_i, \kappa)<f_{\alpha^*_{i+1}}
\restriction [j_i, \kappa)< f_{\alpha^*} \restriction [j_i, \kappa)< f^*
\restriction [j_i, \kappa)$,
\item[$(\gamma)$] for each $i_0<i_1$ we have:\ \ \
$c(\delta_{i_0},\alpha^*_{i_1})<j_{i_1}$, and $c(\alpha^*_{i_0},
\alpha^*_{i_1})<j_{i_1}$, and $c(\alpha^*_{i_1}, \alpha^*)<j_{i_1}$ and
$c(\delta_{i_1},\alpha^*)<j_{i_1}$,  
\item[$(\delta)$] $A_i\subseteq A^*\cap (\alpha^*_i, \delta_i)$,
\item[$(\epsilon)$] $\otp (A_i)=\mu^+_i$,
\item[$(\zeta)$] $A_i\in M_{\delta_i, j_i}$,
\item[$(\eta)$] $B_i\subseteq \lambda_{j_i}$,
\item[$(\theta)$] $\otp (B_i)=\lambda_{j_i}$,
\item[$(\iota)$] $B_\varepsilon\in M_{\alpha^*_i,j_{i}}$ for $\varepsilon<i$,
\item[$(\kappa)$] for every $\alpha\in \bigcup\limits_{\varepsilon\leq i}
A_\varepsilon\cup\{\alpha^*\}$ and $\zeta\leq i$ and $\beta\in B_{\zeta}\cup
\{f^*(j_\zeta)\}$ we have $d(\alpha,\beta)=g^*(j_\zeta)$.
\end{description}
If we succeed then $A=\bigcup\limits_{\varepsilon<\kappa}A_\varepsilon\cup
\{\alpha^*\}$ and $B=\bigcup\limits_{\zeta<\kappa} B_\zeta$ are as
required. During the induction in stage $i$ concerning $(\iota)$ we use just
$\bigvee\limits_{j<\kappa} B_\varepsilon\in M_{\alpha^*_i,j}$. So assume that
the sequence $\langle (j_\varepsilon,A_\varepsilon,B_\varepsilon):\varepsilon
<i\rangle$ has already been defined. 

We can find $j_i(0)<\kappa$ satisfying requirements $(\alpha)$, $(\beta)$,
$(\gamma)$ and $(\iota)$ and such that $\bigwedge\limits_{\varepsilon<i}
\lambda_{j_\varepsilon}<\mu_{j_i(0)}$. Then for each $\varepsilon<i$ we have
$\delta_\varepsilon\in a^{\alpha^*_i}_{j_i(0)}$ and hence
$M_{\delta_\varepsilon,j_\varepsilon}\prec M_{\alpha^*_i,j_i(0)}$ (for
$\varepsilon<i$). But $A_\varepsilon\in M_{\delta_\varepsilon,j_\varepsilon}$
(by clause $(\zeta)$) and $B_\varepsilon\in M_{\alpha^*_i, j_i(0)}$ (for
$\varepsilon<i$), so $\{A_\varepsilon,B_\varepsilon:\varepsilon<i\}\subseteq
M_{\alpha^*_i,j_i(0)}$. Since ${}^{\kappa>}(M_{\alpha^*_i, j_i(0)})\subseteq
M_{\alpha^*_i,j_i(0)}$ (see (ii)),  the sequence $\langle (A_\varepsilon,
B_\varepsilon):\varepsilon<i\rangle$ belongs to $M_{\alpha^*_i, j_i(0)}$.  
We know that for $\gamma_1<\gamma_2$ in $\nacc(C_{\delta_i})$ we have
$c(\gamma_1,\gamma_2)\leq i$ (remember clause (B) and the choice of $c$). As 
$j_i(0)>i$ and so $\mu_{j_i(0)}\geq\mu^+_i$, the sequence 
\[\bar M^*=:\langle M_{\alpha,j_i(0)}:\alpha\in\nacc(C_{\delta_i})\rangle\]
is $\prec$-increasing and $\bar M^* \restriction\alpha\in M_{\alpha,j_i(0)}$
for $\alpha\in\nacc(C_{\delta_i})$ and $M_{\alpha^*_i,j_i(0)}$ appears in
it. Also, as $\delta_i\in \acc(E)$, there is an increasing sequence $\langle
\gamma_\xi:\xi< \mu^+_i\rangle$ of members of $\nacc(C_{\delta_i})$ such that
$\gamma_0=\alpha^*_i$ and $(\gamma_\xi,\gamma_{\xi+1})\cap E\neq\emptyset$,
say $\beta_\xi\in (\gamma_\xi,\gamma_{\xi+1})\cap E$. Each element of
$\nacc(C_\delta)$ is a successor ordinal, so every $\gamma_\xi$ is a successor
ordinal. Each model $M_{\gamma_\xi, j_i(0)}$ is closed under sequences of
length $\leq \mu^+_i$, and hence $\langle\gamma_\zeta:\zeta<\xi\rangle\in
M_{\gamma_\xi, j_i(0)}$ (by choosing the right $\bar C$ and $\delta_i$'s we
could have managed to have $\alpha^*_i=\min (C_{\delta_i})$, $\{\gamma_\xi:
\xi<\mu^+_i\}=\nacc(C_\delta)$, without using this amount of closure).\\
For each $\xi<\mu^+_i$, we know that
\[(H(\chi),{\in}, {<}^*_\chi)\models \mbox{``}(\exists x\in A^*)[
x>\gamma_\xi\ \&\ (\forall \varepsilon<i)(\forall y\in B_\varepsilon)
(d(x,y)=g^*(j_\varepsilon))]\mbox{''}\]
because $x=\alpha^*$ satisfies it. As all the parameters, i.e. $A^*$,
$\gamma_\xi$, $d$, $g^*$ and $\langle B_\varepsilon: \varepsilon<i\rangle$,
belong to $N_{\beta_\xi}$ (remember clauses (e) and (c); note that
$B_\varepsilon\in M_{\alpha^*_i, j_i(0)}$, $\alpha^*_i<\beta_\xi$), there is
an ordinal $\beta^*_\xi\in (\gamma_\xi,\beta_\xi)\subseteq (\gamma_\xi,
\gamma_{\xi+1})$ satisfying the demands on $x$. Now, necessarily for some
$j_i(1,\xi)\in (j_i(0),\kappa)$ we have $\beta^*_\xi\in M_{\gamma_{\xi+1},
j_i(1,\xi)}$. Hence for some $j_i<\kappa$ the set 
\[A_i=:\{\beta^*_\xi: \xi<\mu^+_i\ \ \&\ \ j_i(1,\xi)=j_i\}\]
has cardinality $\mu^+_i$. Clearly $A_i\subseteq A^*$ (as each $\beta^*_\xi\in
A^*$). Now, the sequence $\langle M_{\gamma_\xi, j_i}: \xi<\mu^+_i\rangle
{}^\frown\!\langle M_{\delta_i,j_i}\rangle$ is $\prec$-increasing, and hence
$A_i\subseteq M_{\delta_i, j_i}$. Since $\mu^+_{j_i}>\mu^+_i=|A_i|$ we have
$A_i\in M_{\delta_i,j_i}$. Note that at the moment we know that the set $A_i$
satisfies the demands $(\delta)$--$(\zeta)$. By the choice of $j_i(0)$, as
$j_i> j_i(0)$, clearly $M_{\delta_i, j_i}\prec M_{\alpha^*, j_i}$, and hence
$A_i\in M_{\alpha^*, j_i}$. Similarly, $\langle A_\varepsilon: \varepsilon\leq
i\rangle \in M_{\alpha^*, j_i}$, $\alpha^*\in M_{\alpha^*, j_i}$ and 
\[\sup(M_{\alpha^*, j_i} \cap \lambda_{j_i})=f_{\alpha^*}(j_i)<f^*(j_i).\] 
Consequently,  $\bigcup\limits_{\varepsilon\leq i} A_\varepsilon \cup
\{\alpha^*\}\subseteq M_{\alpha^*, j_i}$ (by the induction hypothesis or
the above) and it belongs to $M_{\alpha^*, j_i}$. Since
$\bigcup\limits_{\varepsilon\leq i} A_\varepsilon \cup \{\alpha^*\}\subseteq 
A^*$, clearly 
\[(H(\chi), \in , <^*_\chi)\models\mbox{``}(\forall x\in\bigcup_{\varepsilon
\leq i} A_\varepsilon \cup \{\alpha^*\})\big(d(x,f^*(j_i))=g^*(j_i)\big)
\mbox{"}.\]
Note that  
\[\bigcup\limits_{\varepsilon\leq i} A_\varepsilon\cup\{\alpha^*\},\ g^*(j_i),
\ d,\ \lambda_{j_i}\in M_{\alpha^*,j_i}\ \mbox{ and }\ f^*(j_i)\in\lambda_{j_i}
\setminus \sup (M_{\alpha^*, j_i}\cap\lambda_{j_i}).\]
Hence the set 
\[B_i=:\{ y<\lambda_{j_i}: (\forall x\in\bigcup_{\varepsilon\leq i} 
A_\varepsilon \cup \{\alpha^*\})(d(x, y)=g^*(j_i))\}\]
has to be unbounded in $\lambda_{j_i}$. It is easy to check that $j_i, A_i,
B_i$ satisfy clauses $(\alpha)$--$(\kappa)$. 

Thus we have carried out the induction step, finishing the proof of the
theorem. \hfill\qed$_{\ref{cztery}}$
\bigskip

\begin{theorem}
\label{piec}
Suppose $\mu$ is singular limit of measurable cardinals. Then
\begin{enumerate}
\item $\left( \begin{array}{c} \mu^+ \\ \mu \end{array} \right)
\rightarrow \left( \begin{array}{c} \mu \\ \mu \end{array}
\right)_\theta$\qquad if $\theta=2$ or at least $\theta<\cf(\mu)$.
\item Moreover, if $\alpha^*<\mu^+$ and $\theta<\cf(\mu)$ then
$\left( \begin{array}{c} \mu^+ \\ \mu \end{array} \right)
\rightarrow \left( \begin{array}{c} \alpha^* \\ \mu \end{array}
\right)_\theta$.
\item If $\theta<\mu$, $\alpha^*< \mu^+$ and $d$ is a function from
$\mu^+\times \mu$ to $\theta$ then for some $A\subseteq \mu^+$,
$\otp(A)=\alpha^*$, and $B=\bigcup\limits_{i< \cf(\mu)} B_i\subseteq
\mu$ we have: 
\[d\restriction A\times B_i\ \mbox{ is constant for each }i<\cf(\mu).\]
\end{enumerate}
\end{theorem}

\Proof Easily 3) $\Rightarrow$ 2) $\Rightarrow$ 1), so we shall prove part 3).

\noindent Let $d: \mu^+\times\mu\rightarrow \theta$. Let $\kappa=:
\cf(\mu)$. Choose sequences $\langle\lambda_i:i<\kappa\rangle$ and $\langle
\mu_i:i<\kappa\rangle$ such that $\langle\mu_i:i<\kappa\rangle$ is increasing
continuous, $\mu=\sum\limits_{i<\kappa} \mu_i$, $\mu_0>\kappa+\theta$, each
$\lambda_i$ is measurable and $\mu_i<\lambda_i<\mu_{i+1}$ (for $i<\kappa$). 
Let $D_i$ be a $\lambda_i$-complete uniform ultrafilter on $\lambda_i$. For
$\alpha< \mu^+$ define $g_\alpha\in {}^\kappa\theta$ by:\quad $g_\alpha(i)=
\gamma$ iff $\{\beta<\lambda_i: d(\alpha,\beta)=\gamma\}\in D$\quad (as
$\theta<\lambda_i$ it exists). The number of such functions is
$\theta^\kappa<\mu$ (as $\mu$ is necessarily strong limit), so for some
$g^*\in {}^\kappa\theta$ the set $A=:\{\alpha<\mu^+:g_\alpha=g^*\}$ is
unbounded in $\mu^+$. For each $i<\kappa$ we define an equivalence relation
$e_i$ on $\mu^+$: 

$\alpha e_i \beta$\quad \underline{iff}\quad $(\forall\gamma<\lambda_i)[
d(\alpha,\gamma)= d(\beta, \gamma)]$.

\noindent So the number of $e_i$-equivalence classes is $\leq
{}^{\lambda_i}\theta<\mu$. Hence we can find $\langle \alpha_\zeta:
\zeta<\mu^+\rangle$ an increasing continuous sequence of ordinals
$<\mu^+$ such that:
\begin{description}
\item[$(*)$] for each $i<\kappa$ and $e_i$-equivalence class X we have:
\item[\qquad\qquad] \underline{either} $X\cap A\subseteq \alpha_0$
\item[\qquad\qquad] \underline{or} for every $\zeta<\mu^+$,
$(\alpha_\zeta, \alpha_{\zeta+1})\cap X \cap A$ has cardinality $\mu$.
\end{description}
Let $\alpha^*=\bigcup\limits_{i<\kappa} a_i$, $|a_i|=\mu_i$, $\langle
a_i: i<\kappa\rangle$ pairwise disjoint. Now we choose by induction on
$i<\kappa$, $A_i$, $B_i$ such that:
\begin{description}
\item[(a)] $A_i\subseteq \bigcup\{(\alpha_\zeta, \alpha_{\zeta+1}):
\zeta\in a_i\}\cap A$ and each $A_i\cap (\alpha_\zeta, \alpha_{\zeta+1})$ is
a singleton,
\item[(b)] $B_i\in D_i$,
\item[(c)] if $\alpha\in A_i$, $\beta\in B_j$, $j\leq i$ then $d(\alpha,
\beta)=g^*(j)$.
\end{description}
Now, in stage i, $\langle (A_\varepsilon, B_\varepsilon): \varepsilon<
i\rangle$ are already chosen. Let us choose $A_\varepsilon$. For each
$\zeta\in a_i$ choose $\beta_\zeta\in (\alpha_\zeta,\alpha_{\zeta+1})\cap A$
such that if $i>0$ then for some $\beta'\in A_0$, $\beta_\zeta e_i \beta'$,
and let $A_i=\{\beta_\zeta: \zeta\in a_i\}$. Now clause (a) is immediate, and
the relevant part of clause (c), i.e. $j<i$, is O.K. Next, as
$\bigcup\limits_{j\leq i} A_j\subseteq A$, the set 
\[B_i=:\bigcap\limits_{j\leq i}\bigcap\limits_{\beta\in A_j}\{\gamma<
\lambda_i: d(\beta, \gamma)=g^*(i)\}\]
is the intersection of $\leq \mu_i<\lambda_i$ sets from $D_i$ and hence
$B_i\in D_i$. Clearly clauses (b) and the remaining part of clause
(c) (i.e. $j=i$) holds. So we can carry the induction and hence finish
the proof. 
\hfill$\qed_{\ref{piec}}$
\bigskip

\bibliographystyle{literal-unsrt}
\bibliography{lista,listb,listx}

\end{document}